\documentclass[12pt]{article}
\usepackage[english]{babel}
\usepackage{amssymb,amsmath,graphics,tikz}

\usepackage{hyperref}

\parskip\medskipamount
\newtheorem{theorem}{Theorem}[section]

\newtheorem{remark}[theorem]{Remark}
\newtheorem{prop}[theorem]{Proposition}
\newtheorem{lemma}[theorem]{Lemma}
\newtheorem{cor}[theorem]{Corollary}
\def\<{\langle}
\def\>{\rangle}
\newcommand{\proof}{\emph{Proof.~}}

\newcommand{\dd}{\mathsf{d}}

\def\qed{{\hfill\hphantom{.}\nobreak\hfill$\Box$}}

\newcommand{\R}{\mathbb{R}}

\numberwithin{equation}{section}

\begin{document}
\title{On metrically complete Bruhat-Tits buildings}
\author{Benjamin Martin \and Jeroen Schillewaert  \and G\"unter F. Steinke  \and Koen Struyve\thanks{This research was supported by a grant from the College of Engineering of the University of Canterbury. The fourth author is supported by  the Fund for Scientific Research ---
Flanders (FWO - Vlaanderen)} }

\maketitle

\begin{abstract}
Metrical completeness for Bruhat-Tits buildings is a natural and useful condition. In this paper we determine which Bruhat-Tits buildings are metrically complete up to certain cases involving infinite-dimensionality and residue characteristic 2.
\end{abstract}



%


%


\section{Introduction}

Non-discrete affine apartment systems, Euclidean buildings or $\R$-buildings were introduced by Jacques Tits (\cite{Tit:86}) as  a non-discrete generalization of affine buildings.  An important subclass of Euclidean buildings are the Bruhat-Tits buildings. One of the main
results of~\cite{Tit:86} shows that any Euclidean building of dimension at least three is a Bruhat-Tits building.
The Bruhat-Tits buildings of dimension at least two have been classified, see~\cite{Tit:86} and \cite{Wei:09}.

Bruhat-Tits buildings are used to study certain groups defined over fields with (not necessarily discrete) valuation.   In general, a Bruhat-Tits building has associated to it an alternative division algebra $A$---often a field
but sometimes a skew field or an octonion algebra---endowed with a valuation.  The building at infinity of a Bruhat-Tits building is a Moufang building.  The root groups of this Moufang building are copies of $A$
or of spaces closely related to $A$, each of which inherits a nonarchimedean metric from the valuation on $A$.  Conversely, given a Moufang building and a nonarchimedean metric
satisfying certain properties on
each of its root groups, there is a unique Bruhat-Tits building which has this Moufang building as its building at infinity.

Euclidean buildings form metric spaces, so the notion of metrical completeness naturally comes into play. While such a property is useful (see for example \cite[Prop.~3.2.4]{Bru-Tit:72} or \cite[Cor.~4.2]{Par:00}), and always satisfied in the discrete case, it does not hold for arbitrary Euclidean buildings.

The aim of the paper is to determine which (non-discrete) Bruhat-Tits buildings of dimension at least two are metrically complete. By a result of Fran\c{c}ois Bruhat and Jacques Tits one can reduce the problem to the study of a nested-ball property for certain metrics defined on the root groups induced by a valuation on the underlying division algebra $A$. In particular, this implies that if the associated Coxeter diagram is simply laced, then metrical completeness of the Bruhat-Tits building is equivalent to the valuation on $A$ being spherically complete, as every root group is a copy of $A$ itself.

We devise two methods to study the other cases: one to study finite-dimensional extension structures and one to study substructures. In Section~\ref{section:conc} we use these tools to discuss the possible two-dimensional cases, which form the building blocks of the higher-dimensional Bruhat-Tits buildings.
We prove the following.



\bigskip
 {\em The Bruhat-Tits building is metrically complete if and only if the completion of the associated alternative division algebra is spherically complete, up to certain cases involving infinite-dimensionality
 and residue characteristic 2.}

\bigskip
A precise version of this statement and a detailed description of each case, as well as a discussion of the exceptions, can be found in Section~\ref{section:conc} (see
Theorems~\ref{thm:forwards} and \ref{thm:cvse}).

\bigskip
\noindent {\em Acknowledgements:} We are grateful to the referees for some corrections and for helpful comments concerning the exposition.


\section{Tools to study metrical completeness}


We make use of a criterion by Bruhat and Tits (\cite[Prop.~7.5.4]{Bru-Tit:72}) which we will explain in detail in Section~\ref{section:rgd} (see Theorem~\ref{thm:compcrit}). In
Sections~\ref{section:xt} and \ref{section:desc} we introduce two tools, one to deal with extensions and a second to deal with descent situations. We apply these tools to the classification in Section~\ref{section:conc}.

\subsection{Root groups of Bruhat-Tits buildings}\label{section:rgd}

Let $X$ be a Bruhat-Tits building and let $Y$ be the building at infinity. For each root $\alpha$ of $Y$, the corresponding root group $U_\alpha$ admits a valuation-like map
$\phi_\alpha: U_\alpha \rightarrow \R \cup \{\infty \} $ obeying the following property (where $e$ is the identity element of $U_\alpha$):
\begin{equation}\label{eq:om}
\forall k \in \R \cup \{\infty \} : U_{\alpha,k} := \phi^{-1}_\alpha([k,\infty])\mbox{ is a subgroup, and } U_{\alpha,\infty} = \{e\}.
\end{equation}
The simplest example that can occur is when $U_\alpha$ is a field (viewed as an additive group) and $\phi_\alpha$ is a (real) complete valuation on $U_\alpha$.

This property implies that the function $\dd: U_\alpha^2 \to \R: (x,y) \mapsto 2^{-\phi_\alpha(x y^{-1})}$ is a metric. Even stronger, it is an \emph{ultrametric}, i.e., it satisfies the \emph{strong triangle inequality}
$$\forall x,y,z \in U_\alpha: \dd(x,z) \leq \max(\dd(x,y),\dd(y,z)).$$

Let us list some basic properties of ultrametric spaces for future reference. For a proof and a detailed overview of these metric spaces we refer to~\cite[\S 1.4]{Kuh:*}.
\begin{lemma}\label{lem:ultra}
Every point inside a ball of an ultrametric space is a center of it; furthermore, if two balls intersect then one contains the other.
\end{lemma}

The balls of the form $\{y\in X \vert \dd(x,y) < r\}$ and $\{y\in X \vert \dd(x,y) \leq r\}$ for $x$ in an ultrametric space $(X,\dd)$ and a positive real number $r$ are both open and closed for the induced topology on $X$. We refer to the first kind as \emph{o-balls} and the second kind as \emph{c-balls}.

Consider the following condition on an ultrametric space $X$:

\bigskip
(MC)\quad Every sequence of nested balls with radius bounded from below by a strictly positive constant has non-empty intersection.
\bigskip

If $X$ is complete then (MC) is equivalent to the following condition:

\bigskip
(MC$'$)\quad Every sequence of nested balls has nonempty intersection.
\bigskip

For the example of a valuation on a field (or more generally an alternative division algebra) one calls the valuation \emph{spherically complete} if the corresponding ultrametric satisfies
(MC$'$).  For fields this is equivalent to the valuation being \emph{maximally complete}, which means that there is no proper extension of the field with valuation having the same residue field and value group (\cite{Kap:42}).

An example of a maximally complete field is the so-called \emph{Hahn-Mal'cev-Neumann series}. Let $k$ be a field and let $G$ be an additive subgroup of $\R$.  Consider the formal power series
of the form $\sum_{i \in I} a_i t^i$ where $I$ is a well-ordered subset of $G$ (well-ordered meaning that every subset has a least element). These power series form a field, the
Hahn-Mal'cev-Neumann series $k((t^G))$. The valuation $\nu$ on $k((t^G))$ is given by
$$ \nu\left(\sum_{i \in I} a_i t^i\right)= {\rm min}\{i\mid a_i\neq 0\}. $$
Note that if we choose $G$ to be the integers we obtain the Laurent series over $k$.

We now give the criterion of Bruhat-Tits (see~\cite[Prop.~7.5.4]{Bru-Tit:72}) for completeness of a Bruhat-Tits building.

\begin{theorem}
\label{thm:compcrit}
 Let $X$ be a Bruhat-Tits building and let $Y$ be the building at infinity.  Then $X$ is metrically complete if and only if for every root $\alpha$ of $Y$, the root group $U_\alpha$ (regarded as an
 ultrametric space) satisfies (MC).
\end{theorem}

Usually one wants to work with the complete building at infinity.  The alternative division algebra and root groups associated to the complete building at infinity are the completions of those associated to the other buildings at infinity. Hence in this case one can replace (MC) with (MC$'$) in the above theorem.

%
%
%
%

\subsection{Extensions and property (MC$'$)}\label{section:xt}

We will study objects of the following form. We start with a skew field $K$ with spherically complete valuation $\nu: K \to \R \cup \{\infty\}$. As explained in the previous section, one can endow this skew field with the structure of a metric space satisfying (MC$'$), which we denote by $(K,\dd_K)$.


Consider an extension $S$ of the additive group of a right vector space $V$ over $K$ by some group $T$. We refer to the group operation of $S$ as addition and we denote it by $\oplus$ although this operation may not be commutative; we also write the inverse of $z\in S$ as $-z$. The identity element of $S$ is denoted by 0.
Let $T'$ be the natural embedding of $T$ in $S$ as a normal subgroup of $S$. 

We also demand that one can lift the scalar product on $V$ to a (right) action $\odot: S \times K \to S$ (so that the action modulo $T'$ is the scalar product on $V$), and moreover that the distributive law with respect to the addition in $S$ holds and that the scalar multiplication is compatible with the field multiplication.
Lastly we ask that the identity element of $K$ acts trivially on $S$.
Explicitly one has the following four equations for the action $\odot$.
\begin{align}
\forall l \in K, x \in S &: ( x \odot l ) \oplus T' = (x \oplus T'). l \label{eq:odot1} \\
\forall l \in K, x, y \in S &: (x \oplus y) \odot l = (x \odot l) \oplus (y \odot l) \label{eq:odot2} \\
\forall l, k \in K, x \in S &: (x \odot k) \odot l = x \odot  (kl) \label{eq:odot3} \\
\forall x \in S &: x \odot 1 =x \label{eq:odot4}
\end{align}

The symbol `.' in the above equations denotes the scalar product on the vector space $V$. Note that $x  \odot (-1)$ is not necessarily the same as $-x$.

Let $\omega$ be a map from $S$ to $\R \cup\{\infty\}$ and $m> 0$ be a real number  satisfying the following three conditions:

\begin{itemize}
\item[(1)] $\forall x \in S: \omega(x)=\infty$ if and only if $x=0$,
\item[(2)] $\forall x \in S, l \in K: \omega( x \odot l)=  \omega(x) + m \nu(l) $, 
\item[(3)] $\forall x, y \in S: \omega(x) = \omega(-x)$ and $\omega(x \oplus y)\geq \min(\omega(x),\omega(y))$.
\end{itemize}

These properties have the following implication.
\begin{lemma}
\label{lem:triangle}
For each $x,y \in S$ one has that $\omega(x) > \omega(y)$ implies $\omega(x \oplus y)= \omega(y)$.
\end{lemma}
\proof
Suppose that both $\omega(x) > \omega(y)$ and  $\omega(x \oplus y) \neq \omega(y)$ hold for certain $x, y \in S$. Condition (3) yields that $\omega(x \oplus y) > \omega(y)$, which in turn implies $\omega(y) = \omega((-x) \oplus x \oplus y) \geq \min(\omega(-x),\omega(x\oplus y)) = \min(\omega(x),\omega(x\oplus y)) > \omega (y)$, a contradiction.
\qed
%

\begin{remark}\label{rem:free} \em
When the group $S$ occurs as a root group $U_\alpha$ of a Bruhat-Tits building and $\omega$ is the valuation-like map $\phi_\alpha$ mentioned at the beginning of Section~\ref{section:rgd}, Conditions (1) and (3) hold automatically, because $\phi_\alpha$ satisfies Equation~\ref{eq:om} from Section~\ref{section:rgd}.  The first part of Equation~\ref{eq:om} is a rewording of Condition (3) while the second part of Equation~\ref{eq:om} is precisely Condition (1).
\end{remark}

%
%

The map $\omega$ induces an ultrametric on $S$ setting $\dd(x,y):=2^{-\omega(x \oplus (- y))}$. The symmetry of $\dd$ follows from the first part of Condition (3) and the strong triangle inequality from Lemma~\ref{lem:triangle}. Note that letting a constant $z\in S$ act on $S$ from the right by addition induces a self-isometry of $S$.

Our goal is to prove property (MC$'$) for $(S, \dd)$ assuming that this property holds for $(T', \dd \vert_{T' \times T'}) $ and that $V$ is finite-dimensional.

If an element $0\neq x \in S$ satisfies the property $$
\omega(x \oplus y)  = \min (\omega(x), \omega(y) ) \mbox{  for all } y \in T'$$
we call it an \emph{independent element}. Note that such an element cannot lie in $T'$. To see this, just take $y$ to be $-x$. 

\begin{lemma}\label{lem:ind_transf}
If $x$ is an independent element then so is $ x \odot l$ for each non-zero $l$ in $K$.
\end{lemma}
\proof
This follows because
\begin{align*}
\omega((x\odot l) \oplus y) &= \omega ((x\odot l) \oplus (y \odot 1)) \\
&= \omega ((x\odot l) \oplus ((y \odot l^{-1})\odot l)) \\
&= \omega ((x \oplus (y \odot l^{-1})) \odot l) \\
&= \omega (x \oplus (y \odot l^{-1})) + m \nu(l) \\
&= \min(\omega (x),\omega(y \odot l^{-1})) + m \nu(l) \\
&= \min(\omega(x \odot l),\omega(y)),
\end{align*}
where the first equality makes use of Equation~\ref{eq:odot4}, the second of~\ref{eq:odot3} and the third of~\ref{eq:odot2}. The fourth equality holds by Condition (2), the fifth because $x$ is an independent element and the sixth by Condition (2) and Equations~\ref{eq:odot3} and \ref{eq:odot4}.
\qed

\begin{lemma}\label{indep-elem}
If $T'$ satisfies property (MC\,$'$) then there exists an independent element in $S$.
\end{lemma}
\proof
Let $x$ be an element of $S \setminus T'$, and let $r:=\inf_{y \in T'} \dd(x,y)$. Since $x$ does not lie in $T'$ and $(T', \dd \vert_{T' \times T'}) $ is metrically complete (as implied by (MC$'$)), it follows that this infimum is non-zero. Let $y_i$ ($i \in \mathbb{N}$) be a sequence of elements of $T'$ such that $\dd(y_i, x)$ converges monotonically to $r$. The c-balls $B_i$ ($i \in \mathbb{N}$) in $T'$ with center $y_i$ and radius $\dd(y_i, x)$ form a sequence of nested balls containing $x$ (by Lemma~\ref{lem:ultra}). Hence, property (MC$'$) for $(T', \dd \vert_{T' \times T'}) $ implies that there exists an element $z \in T'$ in the intersection of all these balls. Clearly, $\dd(x,z) = r$.

We now claim that $x' := x \oplus (-z)$ is an independent element. Note that $r = 2^{-\omega(x \oplus (-z))}=2^{-\omega(x')}$. Suppose that $x'$ is not independent; then there exists $z' \in T'$ such that $\omega(x' \oplus z')  > \min (\omega(x'), \omega(z') )$. By Lemma~\ref{lem:triangle} this is only possible if $\omega(x') =\omega (z' ) < \omega(x' \oplus z')$. 

Calculating the distance between $x$ and $(-z') \oplus z \in T'$ we obtain
\begin{align*}
\dd(x,(-z') \oplus z) &= 2^{-\omega(x \oplus (-z) \oplus z')} \\
&= 2^{-\omega(x' \oplus z')} \\
& <  2^{-\omega(x')} = r.
\end{align*}
This contradicts the way we defined $r$, hence $x'$ is an independent element.
\qed

\begin{lemma}
\label{lem:surjballs}
Assume that $T'$ satisfies property (MC\,$'$) and that $V$ is $K$ considered as a 1-dimensional right vector space over itself. The canonical epimorphism $\rho$ from $S$ to $K$ maps c-balls of $(S,\dd)$ surjectively to c-balls of $(K,\dd_K)$.
\end{lemma}
\proof
The previous lemma allows us to choose an independent element $x$ in $S$. Since $x$ cannot lie in $T'$, one has that $\rho(x) \neq 0$. By 1-dimensionality of $V$ we can express each element $z \in S$ in the  form $z=(x \odot l) \oplus y$ with $l \in K$ and $y \in T'$ in a unique way. Note that $\rho(z) = l \rho(x) $ by Equation~\ref{eq:odot1}.

Fix some element $z$ in $S$ and a c-ball $B$ with center $z$ and radius $r$. Let $z'$ be another point in $B$ and write $z = (x \odot l) \oplus y$ and  $z' = (x \odot l') \oplus y'$ as above. Because $\rho$ is a group homomorphism, $\rho(z) = l \rho(x)$ and $\rho(z') = l' \rho(x)$, we know that $\rho(z \oplus (-z')) = (l-l') \rho(x)$, so we can write $z \oplus (-z')$ as $(x \odot (l  - l' )) \oplus y''$ for some $y'' \in T'$. Hence
\begin{align*}
\dd(z,z') &=2^{-\omega ((x \odot (l - l'))\oplus y'')} \\
&=2^{-\min(\omega (x \odot (l - l')),\omega(y''))} \\
&\geq 2^{-\omega (x \odot (l - l'))} \\
&=2^{ - \omega(x) -m\nu(l-l') }  \\
&= 2^{-\omega(x) +m \nu(\rho(x)) } (\dd_K (\rho(z),\rho(z')))^m,
\end{align*}
where the second equality holds because $ x \odot (l  - l' )$ is again an independent element by Lemma~\ref{lem:ind_transf}.  This implies that $B$ is mapped into a ball $B'$ with center $\rho(z)$ and radius $r' := (2^{\omega(x) -m \nu(\rho(x)) }r)^{1/m}$.

The only thing left to prove is that $\rho(B) = B'$ (i.e., surjectivity). Suppose that $l' \rho(x) \in K$ is an element of the ball $B'$. The element $z' := (x \odot l' ) \oplus y \in S$ is mapped to $l' \rho(x)$ by $\rho$. Moreover, $z \oplus (-z') = x \odot (l  - l' ) $ and hence the inequality above becomes an equality because $y''=0$, yielding $\dd(z,z')= 2^{-\omega(x) +m \nu(\rho(x)) } (\dd_K (\rho(z),\ l' \rho(x)))^m$. In particular this means that $z'$ is an element of the ball  $B$. \qed


\begin{prop} \label{prop:ext}
If the metric space $(T', \dd \vert_{T' \times T'}) $ satisfies (MC\,$'$) and $V$ is finite-dimensional, then the metric space $(S,\dd)$ satisfies (MC\,$'$) as well.
\end{prop}
\proof
We prove this by induction on the dimension of $V$.

The zero-dimensional case is trivial.  The next case is when $V$ is 1-dimensional: here the vector space $V$ over $K$ is $K$ itself as in the statement of Lemma~\ref{lem:surjballs}. Let $B_i$ ($i \in \mathbb{N}$) be a sequence of nested balls in $(S,\dd)$. The images under $\rho$ again form a sequence of nested balls in $(K, \dd_K)$ by the above lemma. As we assumed that $K$ is spherically complete, this sequence contains a common element $l$. Let $z\in S$ be an element such that $\rho(z)=l$. Since adding $-z$ from the right in $S$ is a self-isometry of $(S,\dd)$, this operation produces a new sequence $B'_i$ ($i \in \mathbb{N}$) of nested balls. The images under $\rho$ of this new sequence have to contain the zero element of $K$, or equivalently each ball $B'_i$ contains an element of $T'$.

As the center of a ball in an ultrametric space is arbitrary (see Lemma~\ref{lem:ultra}), one can assume that all the centers lie in $T'$. Because by assumption (MC$'$) holds for $(T', \dd \vert_{T' \times T'}) $, there exists a point in the intersection $\bigcap_{i\in \mathbb{N}} (B'_i \cap T')$, and hence also in $\bigcap_{i\in \mathbb{N}} B_i$.

Now consider the case when $V$ is an $n$-dimensional vector space and suppose that we have proved the result for lower-dimensional vector spaces. Let $U$ be a 1-dimensional subspace of the vector space $V$ and $W$ a complement of $U$ in $V$. Then $V$ is the vector space direct sum of $U$ and $W$. The inverse image $S'$ of $U$ under the canonical epimorphism from $S$ to $V$ is a subgroup of $S$ which is an extension of the additive group of the 1-dimensional right vector space $U$ by the group $T$.

The subgroup $S'$ is closed under the action $\odot$ as it satisfies Equation~\ref{eq:odot1}. The conditions (1)--(3) that are satisfied by $\omega$ are also satisfied by its restriction to $S'$. By the 1-dimensional case above, we have that the metric space $(S', \dd \vert_{S'})$ satisfies (MC$'$).

Interpreting $V$ as the vector space direct sum of $U$ and $W$, one can define a linear projection from $V$ to $W$ with kernel $U$. Combining this with the canonical epimorphism from $S$ to $V$, one obtains a group epimorphism from $S$ to the additive group of $W$ with kernel $S'$. So $S$ can be considered as an extension of the additive group of the right vector space $W$ by the group $S'$. The action $\odot$ satisfies Equations~\ref{eq:odot1}--\ref{eq:odot3} replacing $T'$ with $S'$. As $\omega$ still satisfies Conditions (1)--(3), one can apply the induction hypothesis to the extension of the additive group of $W$ by $S'$ (as $W$ is $(n-1)$-dimensional) and conclude that the metric space $(S,\dd)$ satisfies (MC$'$). This completes the proof.
\qed

Sometimes we will apply the above result in the case when $T$ is the trivial group, reducing to the case when $S$ is a vector space itself. A proof for this special case can also be found in~\cite[Lem.~9.34]{Kuh:*}.

\subsection{Descent and property (MC$'$)}\label{section:desc}

Whereas the last section dealt with proving (MC$'$) for larger structures starting from smaller structures, we now focus on substructures. Consider a vector space $V$ defined over a field $K$ with valuation $\nu$, equipped with a map $\omega$  to $\R \cup \{ \infty \}$ satisfying Conditions (1)--(3) from Section \ref{section:xt} (with $\oplus$ being the usual vector addition, and $\odot$ the scalar product).  We do not need to assume $V$ is finite-dimensional over $K$.

If $U$ and $W$ are non-zero vector subspaces of $V$ such that
\begin{itemize}
\item
$\omega(x +y)  = \min (\omega(x), \omega(y) )$ for all $x\in U$, $y \in W$,
\end{itemize}
we call $U$ and $W$ \emph{$\omega$-complementary}. Note that such subspaces must have zero intersection.
If in addition $V$ is the vector space direct sum of $U$ and $W$ then we call $U$ and $W$ \emph{$\omega$-complements}.

\begin{lemma}\label{sum}
Suppose that the subspaces $U$ and $W$ of $V$ are $\omega$-complements. Then $V$ satisfies (MC\,$'$)
if and only if $U$ and $W$  satisfy (MC\,$'$).

\end{lemma}
\proof
Write $x=u+w$ and $x'=u'+w'$ with $u,u'\in U$, $w,w'\in W$. Then, because $U$ and $W$ are $\omega$-complements, one has
$$\dd(x,x')=\max(\dd(u,u'),\dd(w,w')).$$
Hence $V$ is isometric to the direct product of the ultrametric spaces  $U$ and $W$. This means we can decompose any ball in $V$ into a Cartesian product of a ball in $U$ with a ball in $W$. The statement then follows.
\qed

%
%
Let $F$ be a subfield of $K$ and let $\sigma:V \to V $ be an $F$-linear map preserving $\omega$ such that $\sigma^2={\rm id}_V$, and let $V_\sigma := \{v\in V \vert v^\sigma = v \}$ be the vector subspace of $V$ it fixes.
Suppose that the characteristic of the residue field of $(K,\nu)$ is not 2 (so that $\nu(2) = 0$). Denote by $V'_\sigma$ the set $\{v\in V \vert v^\sigma = -v \}$. The next lemma allows us to decompose the vector space $V$.

\begin{lemma}\label{V_sigma}
$V_\sigma$ and $V'_\sigma$ are $\omega$-complements in $V$.
\end{lemma}
\proof
For $x\in V$ the vectors $v=\frac12(x+x^\sigma)$ and $v'=\frac12(x-x^\sigma)$ belong to $V_\sigma$ and $V'_\sigma$, respectively. Observe that $x=v+v'$. As $V_\sigma \cap V'_\sigma = \{0\}$ this implies that $V$ is the vector space direct sum of $V_\sigma$ and $V'_\sigma$. Furthermore
\begin{eqnarray*}
\omega(x)&\geq&\min(\omega(v),\omega(v'))\\
&=&\min(-m\nu(2)+\omega(x+x^\sigma),-m\nu(2)+\omega(x-x^\sigma))\\
&\geq&\min(\min(\omega(x),\omega(x^\sigma)),\min(\omega(x),\omega(-x^\sigma)))\\
&=&\min(\omega(x),\omega(x^\sigma))\\
&=&\omega(x).
\end{eqnarray*}
Hence $\omega(v+v')=\min(\omega(v),\omega(v'))$ and $V_\sigma$ and $V'_\sigma$ are $\omega$-complementary.
\qed

\begin{cor}\label{cor:des}
If $V$ satisfies (MC\,$'$) then so does $V_\sigma$.
\end{cor}
\proof
This follows immediately from Lemmas~\ref{sum} and~\ref{V_sigma}.
\qed

\section{Results}\label{section:conc}
Recall from Theorem~\ref{thm:compcrit} that a Bruhat-Tits building is metrically complete if and only if property (MC) holds for the metric space associated to each root group. Using the classification of Bruhat-Tits buildings of dimension at least 2 (see~\cite{Tit:86},~\cite{Wei:09} and~\cite{Hit-Kra-Wei:10}), one observes that there are at most two different types of root group in the building at infinity for a given Bruhat-Tits building. Accordingly, there are at most two isometry classes of associated metric spaces we need to consider. Moreover, it follows from the classification that there is a rank 2 residue where the different types of root group occur simultaneously. Such a residue we will call a \emph{highest label residue} (as these residues correspond with the edges of the Coxeter diagram with the highest label). All highest label residues have the same isomorphism type. One can look up this isomorphism type in Sections 12.13 up to 12.19 of~\cite{Wei:03}.

This observation allows us to reduce our discussion to the Bruhat-Tits buildings of dimension 2.
The buildings at infinity here are generalized Moufang polygons, of which a detailed description can be found in~\cite[\S 16]{Tit-Wei:02}. There are nine types of Moufang polygons, which in the notation of \cite{Tit-Wei:02} are
$\mathcal{T}(A)$,
$\mathcal{Q_I}(K, K_0, \sigma)$,
$\mathcal{Q_Q}(K, L_0, q)$,
$\mathcal{Q_D}(K, K_0, L_0)$,
$\mathcal{Q_P}(K, K_0, \sigma, L_0,q)$,
$\mathcal{Q_E}(K, L_0, q)$,
$\mathcal{Q_F}(K, L_0, q)$,
$\mathcal{H}(J,F, \#)$, and
$\mathcal{O}(K, \sigma)$.
One can associate an alternative division algebra to each of these: namely, whichever one of  $A$, $K$ and $F$ is applicable. For the Bruhat-Tits buildings associated to these we refer to~\cite[\S 19-25]{Wei:09} and~\cite{Hit-Kra-Wei:10}.


Let us now give a precise version of the result stated in the introduction. For clarity of the statements we consider the complete building at infinity. The direction going from the Bruhat-Tits building to the associated alternative division algebra can be stated as follows.


\begin{theorem}\label{thm:forwards}
Let $X$ be a Bruhat-Tits building of dimension at least 2 and let $Y$ be a highest label residue of the complete building at infinity.  Suppose $X$ is metrically complete.
\begin{enumerate}
\item
If $Y$ is one of the generalized Moufang polygons
\begin{itemize}
\item
$\mathcal{T}(A)$,
\item
$\mathcal{Q_I}(K, K_0, \sigma)$,
\item
$\mathcal{Q_Q}(K, L_0, q)$,
\item
$\mathcal{Q_P}(K, K_0, \sigma, L_0,q)$,
\item
$\mathcal{H}(J,F, \#)$, or
\item
$\mathcal{O}(K, \sigma)$
\end{itemize}
then $A$, $K$ or $F$ (where applicable) is spherically complete.
\item
If $Y$ is the generalized Moufang polygon $\mathcal{Q_E}(K, L_0, q)$ then the associated field $K$ is spherically complete provided that the
residue characteristic of $K$ is different from 2.
\item
If $Y$ is one of the generalized Moufang polygons $\mathcal{Q_D}(K, K_0, L_0)$ or
$\mathcal{Q_F}(K, L_0, q)$ then no conclusion can be drawn on the spherical completeness of the associated field $K$.
\end{enumerate}
\end{theorem}


%
%
The other direction can be stated as follows.

%

\begin{theorem}\label{thm:cvse}
Let $X$ be a Bruhat-Tits building of dimension at least 2 with $Y$ a highest label residue of the building at infinity. Suppose that the alternative division algebra associated to $Y$ is spherically complete.
Then $X$ is metrically complete if $Y$ is the generalized Moufang polygon \dots under the assumption that \dots .  (Here the dots indicate the appropriate entry from the table below,
and a `/' for the assumption means that there is no extra assumption.)

\begin{center}\begin{tabular}{c|c} Generalized polygon & assumptions  \\\hline
 $\mathcal{T}(A)$ & / \\
$\mathcal{Q_I}(K, K_0, \sigma)$ & residue characteristic different from 2 \\
$\mathcal{Q_Q}(K, L_0, q)$ & $L_0$ finite-dimensional \\
$\mathcal{Q_D}(K, K_0, L_0)$ & $K$ finite-dimensional over $K^2$ \\
$\mathcal{Q_P}(K, K_0, \sigma, L_0,q)$ & residue characteristic different from 2, \\
 & $L_0$ finite-dimensional \\
$\mathcal{Q_E}(K, L_0, q)$ &  / \\
$\mathcal{Q_F}(K, L_0, q)$& $K$ finite-dimensional over $K^2$\\
$\mathcal{H}(J,F, \#)$&  $J$ finite-dimensional\\
$\mathcal{O}(K, \sigma)$&  /
\end{tabular}
\end{center}
\end{theorem}

%

Concerning the extra assumptions for both directions: the restriction to finite-dimensional cases is sometimes needed as will be illustrated in the quadratic form type case
$\mathcal{Q_Q}(K, L_0, q)$ below. The necessity of having residue characteristic different from 2 where needed below is unclear to the authors. Some special cases will be discussed for the involutory type case $\mathcal{Q_I}(K, K_0, \sigma)$.

In the following we discuss each case separately and in detail. Again note that the building at infinity is implied to be the complete building at infinity. 

\paragraph{The triangles $\mathcal{T}(A)$.}
This case corresponds with the simply laced case mentioned in the introduction.
These are parametrized by an alternative division algebra $A$ with a fixed valuation.  Each root group is the additive group of $A$ with the given valuation. So we obtain directly from
Theorem~\ref{thm:compcrit} that the  Bruhat-Tits building is metrically complete if and only if this valuation is spherically complete.

\paragraph{The quadrangles $\mathcal{Q_I}(K, K_0, \sigma)$ of involutory type.}

These are parametrized by a skew field $K$ equipped with a valuation and the fixed subset $K_0 \subset K$ of a valuation-preserving involution $\sigma$ (here we assume the characteristic of $K$ is not 2; in the characteristic 2 case one has more choices for $K_0$).  Two types of root group occur for each such quadrangle: the additive group of $K$ and the additive group of $K_0$.  

If the Bruhat-Tits building is metrically complete, one obtains directly from Theorem~\ref{thm:compcrit} that the valuation on $K$ is spherically complete. On the other hand, if the
valuation on $K$ is spherically complete and has residue characteristic different from 2, one can apply Corollary~\ref{cor:des} taking the subfield $F$ to be the fixed field of $\sigma$ in $Z(K)$ to obtain that the metric space on $K_0$ satisfies (MC$'$). Then one
 deduces from Theorem~\ref{thm:compcrit} that the Bruhat-Tits building is metrically complete.

Let us consider the special case when $K$ is a field and $K_0$ is the subset of $K$ fixed by $\sigma$. In particular this implies that $K_0$ is a subfield. An open question in this special case is whether a maximally complete field $K$ with residue characteristic 2 can admit an involution-preserving valuation such that the fixed subfield is no longer maximally complete. In two cases it is possible to exclude this.


First assume that $K$ is separably closed and not of characteristic 2. Since $K$ is a quadratic extension of $K_0$, the results of Artin and Schreier (\cite{Art-Sch:26}, \cite{Art-Sch:27}) imply that the characteristic of $K$ is zero, and that $K_0$ is an ordered field. Hence $K_0$ contains $\mathbb{Q}$, and since $K_0$ is complete (it is the fixed point set of a valuation-preserving involution of a complete field), it also contains a copy of the $2$-adic numbers. However, the field $\mathbb{Q}_2$ cannot be ordered (as it contains a square root of $-7$), hence we have a contradiction.

Note that the results of Artin and Schreier also allow one to construct an involution for each maximally complete algebraically closed field of characteristic zero. Nevertheless this involution does not preserve the valuation and the fixed subfield is not maximally complete---or even metrically complete---for the reasons outlined above.

Another case where such involutions can be excluded is when $K$ is a Hahn-Mal'cev-Neumann series $k((t^\mathbb{Q}))$ over a perfect field $k$ with characteristic 2 having an extension of degree divisible by 2 (these fields are defined at the end of Section~\ref{section:rgd}). This follows from results of Kedlaya and Poonen (\cite{Ked-Poo:05}), which show among other things that the fixed subfield is a Hahn-Mal'cev-Neumann series over a subfield of $k$, and hence again maximally complete.

\paragraph{The quadrangles $\mathcal{Q_Q}(K, L_0, q)$ of quadratic form type.}

These are parametrized by a field $K$ with a valuation $\nu$ and a vector space $L_0$ over $K$ with an anisotropic quadratic form $q$.  The valuation-like map on $L_0$ is defined by $\nu \circ q$.

Two types of root group occur for each such quadrangle: the additive groups of $K$ and $L_0$ respectively with the given valuations.  One obtains directly from Theorem~\ref{thm:compcrit}
that if the Bruhat-Tits building is metrically complete, then the field is maximally complete. The other direction is also true by Theorem~\ref{thm:compcrit} provided that $L_0$ is finite-dimensional.

To show this we use the results of Section~\ref{section:xt}, where $S$ is the additive group of the vector space $L_0$, the group $T$ is trivial and the map $\omega$ is given by $\nu \circ q$. Conditions (1) and (2) (with $m=2$) from Section~\ref{section:xt} are satisfied as $q$ is an anisotropic quadratic form. Since the Moufang polygon is the building at infinity of a
Bruhat-Tits building, it can be shown that Condition (3) is satisfied (see ~\cite[Prop. 19.4; Thm. 19.18]{Wei:09} and Remark~\ref{rem:free}). Hence we can apply Proposition~\ref{prop:ext}
to $L_0$.

Let us give a counterexample in the infinite-dimensional case. Consider as field $K$ the Hahn-Mal'cev-Neumann series $\R((t^\R))$ with the usual valuation $\nu$. This field is maximally complete. Let $V$ be an infinite-dimensional vector space over $K$ with basis $e_0, e_1, e_2, \dots .$ The map $q : V \to K: (a_0, a_1, a_2, \dots) \mapsto a_0^2 + a_1^2 + a_2^2 +\dots $ (with respect to this basis) is an anisotropic quadratic form. It it easy to show that $\nu(q(a_0,a_1,\dots )) = 2 \min_i (\nu(a_i))$.  This implies that $\nu \circ q$ satisfies Condition (3) from Section~\ref{section:xt} (as mentioned above, Conditions (1) and (2) are satisfied automatically). In particular one can associate a Bruhat-Tits building to the quadrangle $\mathcal{Q_Q}(K, V, q)$.  

We now claim that this building is not metrically complete. Using the explicit formula for $\nu \circ q$ one can express the distance function $\dd((a_0,a_1,\dots),(b_0,b_1, \dots))$  defined in Section~\ref{section:xt} as $\max_i (4^{-\nu(a_i -b_i)})$. Consider the sequences $x_i$ and $r_i$ where
$$
\begin{array}{ll}
x_0 := (1, 0, 0, 0, \dots) &  r_0 = 4^{-1/2} \\
x_1 := (1, t^{1/2}, 0, 0,  \dots) &  r_1 = 4^{-3/4}\\
x_2 := (1, t^{1/2} , t^{3/4}, 0, \dots) & r_2 = 4^{-7/8} \\
\dots &
\end{array}
$$
The balls $B_i$ with center $x_i$ and radius $r_i$ contain all vectors $x_{i'}$ with $i' \geq i$, hence these balls are nested by the basic properties of ultrametric spaces mentioned in Lemma~\ref{lem:ultra}. Suppose that all of these balls have a vector $c :=(c_0, c_1, \dots)$ in common. As only a finite number of the $c_i$ are non-zero, one can let $j$ be the first coordinate such that $c_j$ is zero. This implies that $\dd(x_j,c)$ is at least $4^{-(1-1/2^j)}$, so $c \notin B_j$, which contradicts the way we constructed $c$.  Hence $V$ does not satisfy property (MC). It follows from Theorem~\ref{thm:compcrit} that the Bruhat-Tits building is not metrically complete.

\paragraph{The quadrangles $\mathcal{Q_D}(K, K_0, L_0)$ of indifferent type.}

These are parametrized by certain subsets $K_0$ and $L_0$ of a field $K$ of characteristic 2 with valuation; $K_0$ and $L_0$ can be viewed as vector spaces over the field of squares $K^2$.  Moreover $K_0^2$ generates $K$ as a ring. Two types of root group occur for each such quadrangle: the additive groups of $K_0$ and $L_0$ respectively, equipped with the restriction of the valuation on $K$. Theorem~\ref{thm:compcrit} and Proposition~\ref{prop:ext} imply that the Bruhat-Tits building will be metrically complete if $K$ is maximally complete and $K_0$ and $L_0$ are both finite-dimensional over $K^2$.  (Note that the field isomorphism $x\mapsto x^2$ from $K$ to $K^2$ takes c-balls onto c-balls, so $K^2$ is spherically complete if and only if $K$ is.)

If $K_0$ is finite-dimensional over $K^2$ then $K$ is an extension of $K^2$ generated by a finite number of algebraic elements. Hence $K$ is then finite-dimensional over $K^2$. Conversely if $K$ is finite-dimensional over $K^2$, then the subspaces $K_0$ and $L_0$ are both finite-dimensional over $K^2$. So the Bruhat-Tits building will be metrically complete if $K$ is maximally complete and finite-dimensional over $K^2$.

\paragraph{The quadrangles $\mathcal{Q_P}(K, K_0, \sigma, L_0,q)$ of pseudo-quadratic form type.}

The situation here is similar to the quadrangles of involutory type: $K$, $\sigma$ and $K_0$ are as before. One type of root group is again the additive group of a skew field $K$ with valuation $\nu$. The other type of root group is not $K_0$ but an extension $S$ of a vector space $L_0$ over $K$ by $K_0$.

Let us describe this extension in more detail (using the notation from the beginning of Section~\ref{section:xt}). The elements of the group $S$ are given by $\{(u,t) \in L_0 \times K \vert q(u) - t \in K_0\}$ and the group law is given by $(u,t) \oplus (v,s) = (u+v, t+s+f(v,u))$, where $q$ is an anisotropic pseudo-quadratic form on $L_0$ and $f$ the skew-Hermitian form associated to it (for a definition see~\cite[Def. 11.15]{Tit-Wei:02}). The normal  subgroup $T' := \{(0,t)\mid t\in L_0\}$ of $S$ is isomorphic to $K_0$ and $S/T'$ is canonically isomorphic to $L_0$. So we are indeed dealing with an extension of the vector space $L_0$ over $K$ by $K_0$. The valuation-like map on $S$ is defined by $\omega((u,t)) = \nu(t)$. The right action of $K$ on $S$ can be defined as $ (u,t)\odot k = (u k, k^\sigma t k)$. Using the definition of pseudo-quadratic forms one shows that this product satisfies Equations (2.2)--(2.5). The anisotropy of $q$
means that for any $u \in L_0$, one has that $q(u) \in K_0$ if and only if $u=0$. Conditions (1) and (3) from Section~\ref{section:xt} are automatically satisfied by Remark~\ref{rem:free}.
As in the involutory type case, the involution $\sigma$ preserves the valuation $\nu$. From this one derives that Condition (2) holds with $m=2$.

%

Hence we may conclude by Theorem~\ref{thm:compcrit} that if the Bruhat-Tits building is metrically complete, then the skew field $K$ is spherically complete. The other direction is also true
under the assumptions of residue characteristic different from 2 and finite-dimensionality of $L_0$---for the metric space $K_0$ satisfies (MC$'$) if the residue characteristic is different from 2 by the same reasoning as for the quadrangles of involutory type.
If $L_0$ is finite-dimensional then $S$ satisfies (MC$'$), by Proposition~\ref{prop:ext}.

\paragraph{The quadrangles $\mathcal{Q_E}(K, L_0, q)$ of exceptional type $\mathsf{E}_6$, $\mathsf{E}_7$ and $\mathsf{E}_8$.}

These are parametrized by a finite-dimensional vector space $L_0$ over a field $K$ with valuation $\nu$, using a special kind of anisotropic quadratic form $q$ on $L_0$ to obtain the
valuation-like map $\omega$ on $L_0$.  Two types of root group occur for each such quadrangle: the first is the additive group of $L_0$, as in the quadratic form case, and the second is an extension $S$ of a finite-dimensional vector space $X_0$ over $K$ by $K$. The space $X_0$ is equipped with a binary map $g: X_0 \times X_0 \to K$ and a map $\pi: X_0 \to L_0$.

More precisely this extension is defined on the set $X_0 \times K$ with group law $(a,s) \oplus (b,t) =(a+b, s+t+g(a,b))$. The valuation-like map $\omega$ as in Section~\ref{section:xt} on this extension is given by $(a,t) \mapsto \nu(q(\pi(a)) +t)$, and one can define an action on this extension by setting $(a,t) \odot s := (sa,s^2 t)$. This action satisfies Condition (2) from Section~\ref{section:xt} with $m=4$
by~\cite[Prop 21.24]{Wei:09}. Conditions (1) and (3) are automatically satisfied: see Remark~\ref{rem:free}.

If the field $K$ is maximally complete, we obtain from Theorem~\ref{thm:compcrit} that the Bruhat-Tits building is metrically complete by applying Proposition~\ref{prop:ext} twice: once to
$L_0$ (as in the quadratic form case) and once to $S$.

If the residue characteristic is different from 2, then one can find a basis of $L_0$  such that $q$ becomes of the form $q(x_1, \dots, x_n) = a_1 x_1^2 + \dots a_n x_n^2$ (see for instance~\cite[42:1]{O'm:73}). The involutory map $\phi:(x_1, x_2, \dots, x_n) \mapsto (x_1, -x_2,\dots, -x_n)$ leaves the quadratic form invariant. As fixed set we obtain a 1-dimensional vector space over $K$, which we can identify with $K$ itself. So, using Theorem~\ref{thm:compcrit} and Corollary~\ref{cor:des}, we obtain that
$K$ is maximally complete if the Bruhat-Tits building is metrically complete.

\paragraph{The quadrangles $\mathcal{Q_F}(K, L_0, q)$ of exceptional type $\mathsf{F}_4$.}
As in the indifferent type case,
we work with a field $K$ of characteristic 2. The other parameters are a vector space $L_0$ equipped with an anisotropic quadratic form with certain special properties. Two types of root group occur.  The first is the additive group of $F^4 \times K$ considered as vector space over $K$, where $F$ is a subfield of $K$ containing the field of squares $K^2$. As one can identify $K^2$ with $K$, one can indeed consider $F$ as a vector space over $K$. The second is
the additive group of $K^4 \times F$ considered as vector space over $F$.  Each is equipped with an anisotropic quadratic form. (The fourth powers in the above description denote direct products of four copies of the vector space, not the subfield of fourth powers.)

If these vector spaces are finite-dimensional and the field $K$ is maximally complete, then the Bruhat-Tits building is metrically complete by Theorem~\ref{thm:compcrit} and Proposition~\ref{prop:ext}, as in the quadratic form case. These vector spaces are finite-dimensional if and only if $K$ is finite-dimensional over $F$ and $F$ finite-dimensional over $K^2$. This happens exactly when $K$ is finite-dimensional as vector space over $K^2$.


\paragraph{The hexagons $\mathcal{H}(J,F, \#)$.}
Starting from a quadratic Jordan algebra $J$ of degree 3 over a field $F$, one constructs two types
of root group: the first is the additive group of the field $F$ with valuation and the second is the additive group of $J$ equipped with a valuation-like map $\omega := \nu \circ N$, where $N$ is a certain cubic norm $N$ (so $N(ta)=t^3N(a)$). Except in one case these vector spaces are finite-dimensional. So in all other cases the Bruhat-Tits building is metrically complete if and only if the field $F$ is maximally complete, by Theorem~\ref{thm:compcrit} and Proposition~\ref{prop:ext} (Conditions (1) and (3) for
$\omega$ are again satisfied by Remark~\ref{rem:free}, Condition (2) with $m=3$).

The one case where infinite-dimensionality may occur is the case of a hexagonal system $(E/F)^0$ of type $1/F$, where $\mathsf{char}(F)=3$ and $E$ is a purely inseparable extension of $F$ (see~\cite[(15.20)]{Tit-Wei:02}). Here the Bruhat-Tits building is metrically complete if and only if the fields $F$ and $E$ are both maximally complete, by Theorem~\ref{thm:compcrit}. This
last condition is automatically fulfilled if $F$ is a maximally complete field and the extension is of finite index, by Proposition~\ref{prop:ext}.

\paragraph{The octagons $\mathcal{O}(K, \sigma)$.}
The ingredients here are a field $K$ and a Tits endomorphism $\sigma$ (so $(x^\sigma)^\sigma = x^2$). In order for the octagon to be the building at infinity of a Bruhat-Tits building, the
valuation $\nu$ on $K$ must satisfy $\nu(x^\sigma) = \sqrt{2} \nu(x)$ (see~\cite{Hit-Kra-Wei:10}).

Again there are two types of root group in each case: the first is the additive group of a field $K$ with valuation $\nu$ and the second is the additive group of an extension $S$ of the additive group of $K$ by itself.

This extension is defined on the set $K \times K$ and has group law $(t,u) \oplus (s,v) = (t+s+u^\sigma v, u+v)$. The valuation-like map $\omega$ is defined on $S$ by setting $\omega((t,u)) := \nu(t^{\sigma+2} + ut +u^\sigma)$. We define an action of $K$ on $S$ by $(t,u)\odot k = (tk^{\sigma-1}, uk)$, which satisfies Condition (2) from Section~\ref{section:xt} with $m = \sqrt{2}$. Conditions (1) and (3) hold by Remark~\ref{rem:free}, so one can apply Proposition~\ref{prop:ext} to this extension.

One concludes from Theorem~\ref{thm:compcrit} that the Bruhat-Tits building is metrically complete if and only if the field $K$ is maximally complete.

\bigskip
This completes the proofs of Theorems~\ref{thm:forwards} and \ref{thm:cvse}.

\newpage

B. Martin, Department of Mathematics, University of Auckland, Private Bag 92019, Auckland 1142, New Zealand.  Email: {\em Ben.Martin@auckland.ac.nz} \\ \smallskip

J. Schillewaert, Department of Mathematics, Universit\'e Libre de Bruxelles, CP 216, Bd du Triomphe, 1050 Brussels, Belgium.  Email: {\em jschille@ulb.ac.be}  \\ \smallskip

G. F. Steinke, Department of Mathematics and Statistics, University of Canterbury, Private Bag 4800, Christchurch 8140, New Zealand.  Email: {\em Gunter.Steinke@canterbury.ac.nz} \\ \smallskip

K. Struyve, Department of Mathematics, Ghent University, Krijgslaan 281, 9000 Ghent, Belgium.  Email: {\em kstruyve@cage.ugent.be} \\

\end{document}